\documentclass[12pt]{article}
\usepackage{amsmath,amssymb,amsthm,amsfonts}
\ifx\pdftexversion\undefined
  \usepackage[dvips]{graphicx}
\else
  \usepackage[pdftex]{graphicx}
\fi
\swapnumbers
\usepackage[numeric]{amsrefs}
\usepackage{hyperref}

\newtheorem{thm}[equation]{Theorem}
\newtheorem{cor}[equation]{Corollary}
\newtheorem{lem}[equation]{Lemma}

\newtheorem{definition}[equation]{Definition}

\newcommand{\bx}{\hfill$\square$\vspace{.6cm}}

\DeclareMathOperator{\sg}{sgn}

\numberwithin{equation}{section}

\newcommand{\gobble}[1]{}
  \newcommand{\rangeref}[2]{%
    \ref{#1}--\afterassignment\gobble\fam 0\ref{#2}%
  }

\renewcommand\a{\alpha}         
\renewcommand\b{\beta}
\newcommand\g{\gamma}
\renewcommand\d{\delta}
\newcommand\e{\epsilon}
\renewcommand\l{\lambda}

\newcommand\G{\Gamma}

\newcommand{\N}{{\mathbb{N}}}
\newcommand{\Z}{{\mathbb{Z}}}
\newcommand{\R}{{\mathbb{R}}}

\newcommand{\C}{{\mathbb{C}}}

\newcommand{\Q}{{\mathbb{Q}}}

\newcommand\re{\mbox{Re~}}

\renewcommand\Re{\operatorname{Re}}
\renewcommand\Im{\operatorname{Im}}

\begin{document}

\title{The Highly Oscillatory Behavior \\ of Automorphic Distributions for $SL(2)$}
\author{Stephen D. Miller\thanks{Supported by NSF grant DMS-0301172 and an Alfred P. Sloan Foundation Fellowship}~~and Wilfried
Schmid\thanks{Supported in part by NSF grant DMS-0070714}}

\maketitle

\begin{abstract}
Automorphic distributions for $SL(2)$ arise as boundary values of modular forms and, in a more subtle manner, from Maass forms. In the case of modular forms of weight one or of Maass forms, the automorphic distributions have continuous first antiderivatives. We recall earlier results of one of us on the H\"older continuity of these continuous functions and relate them to results of other authors; this involves a generalization of classical theorems on Fourier series by S. Bernstein and Hardy-Littlewood. We then show that the antiderivatives are non-differentiable at all irrational points, as well as all, or in certain cases, some rational points. We include graphs of several of these functions, which clearly display a high degree of oscillation. Our investigations are motivated in part by properties of ``Riemann's nondifferentiable function", also known as ``Weierstrass' function".
\end{abstract}

\section{Introduction}\label{introduction}

Riemann is credited -- inaccurately perhaps -- with providing the first example of a continuous function which fails to be differentiable at ``most" points:
\begin{equation}
\label{Riemann}
f(x)\ = \ {\sum}_{n\geq 1}\ \frac 1{n^2}\,\sin(2\pi n^2x)
\end{equation}
is non-differentiable except at points $x=p/2q$ with $p$ and $q$
odd; at those, the derivative exists and is equal to $-\pi$. Many
authors have studied this function, beginning with Hardy
\cite{Har}; the final detail was put into place only in 1971
\cite{Ge}. Duistermaat \cite{Du} recounts this literature. He also
gives new proofs of the main properties of this function. His
starting point is the observation that $f'(x)$ exists as a
distribution which is automorphic in an appropriate sense.

The function (\ref{Riemann}) is merely the tip of an iceberg. In this note, we continue the study, begun in \cite{S1}, of the properties of automorphic distributions for subgroups of finite index $\G \subset SL(2,\Z)$. These automorphic distributions have continuous anti-derivatives which are {\em non-differentiable everywhere}, or {\em everywhere with the exception certain rational points}, as in the case of the function (\ref{Riemann}). We establish more: the continuous antiderivatives satisfy global H\"older conditions $|f(y)-f(x)|= O(|y-x|^\alpha)$, but definitely violate the pointwise H\"older conditions $|f(x)-f(x_0)|= O(|x-x_0|^\gamma)$, $\beta<\gamma\leq 1$, for values $\beta=\beta(x_0)\geq\alpha$ which depend on the arithmetic properties of $x_0$. This behavior reflects a high degree of oscillation around all rational points.
\begin{figure}
\begin{center}
\includegraphics[height=2in,width=5.3in]{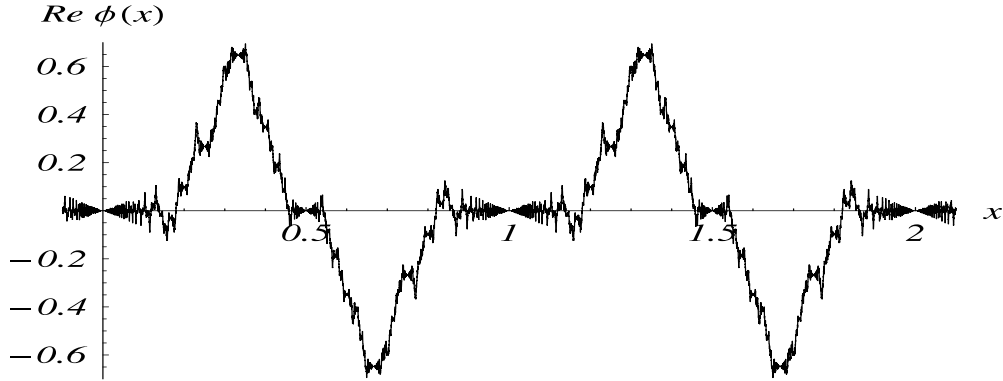}
\end{center}
\caption{The real part of the antiderivative $\phi(x)$ of the automorphic distribution corresponding to the Maass form for $SL(2,\Z)$ with $\lambda\approx 27.56\, i$.}
\label{realmaassbig+_3}
\end{figure}
Figure \ref{realmaassbig+_3}, for example, plots the real part of the antiderivative $\phi(x)$ of the automorphic distribution corresponding to the Maass form of smallest non-zero eigenvalue for $\G=SL(2,\Z)$; $\Re \phi(x)$ is continuous, but everywhere non-differentiable. Near the origin $\phi(x) \sim |x|^{1+\l} \phi (1/x)$, with $\l\approx 27.56\,i$, and this behavior is replicated at all rational points. The absolute value of $\phi(x)$ also oscillates rapidly, as is illustrated by figure \ref{absmaassbig+_3}.
\begin{figure}
\begin{center}
\includegraphics[height=2in,width=5.3in]{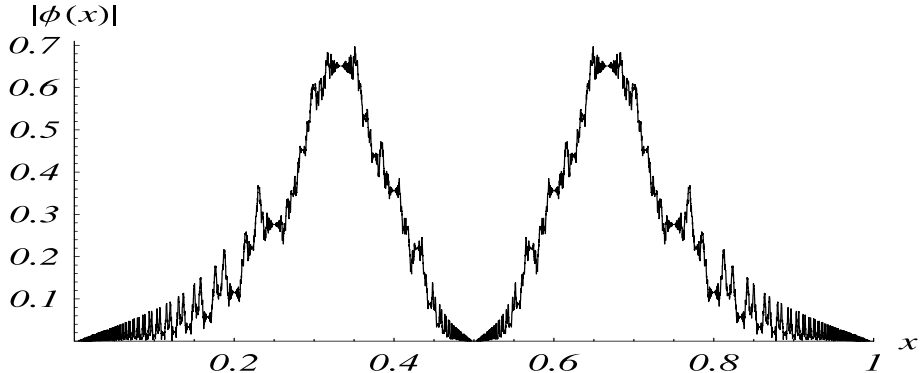}
\end{center}
\caption{The absolute value of of same function $\phi(x)$ as in figure \ref{realmaassbig+_3}.}
\label{absmaassbig+_3}
\end{figure}
Near the origin $|\phi(x)|$ evidently displays fractal behavior~-- see figure \ref{absmaasssmall+_3}.
\begin{figure}
\begin{center}
\includegraphics[height=2in,width=5.3in]{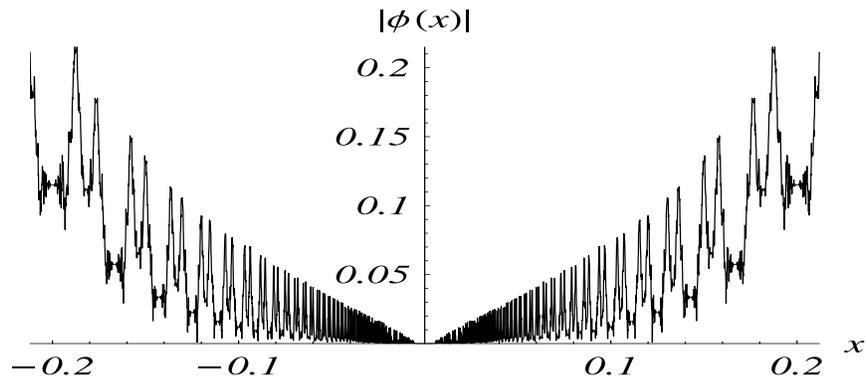}
\end{center}
\caption{Same as figure \ref{absmaassbig+_3}, but near the origin.} \label{absmaasssmall+_3}
\end{figure}

\begin{figure}
\begin{center}
\includegraphics[height=2in,width=5.3in]{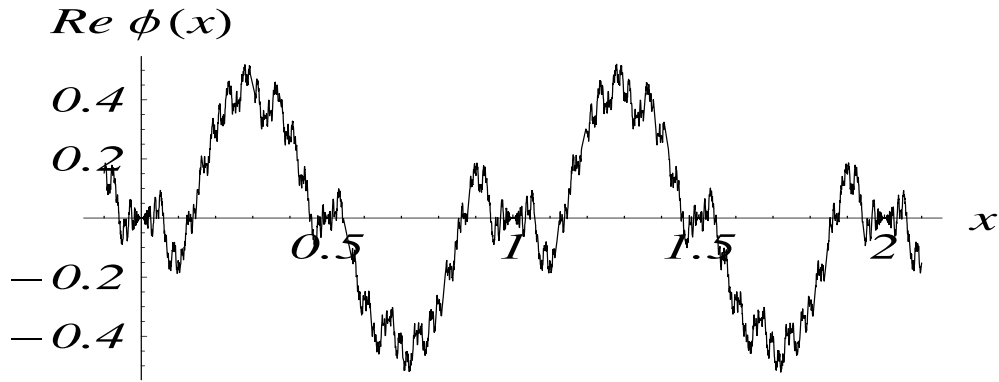}
\end{center}
\caption{The real part of the antiderivative of the boundary distribution of the weight one modular form (\ref{weight1}).} \label{weight1big+_3}
\end{figure}

\begin{figure}
\begin{center}
\includegraphics[height=2in,width=5.3in]{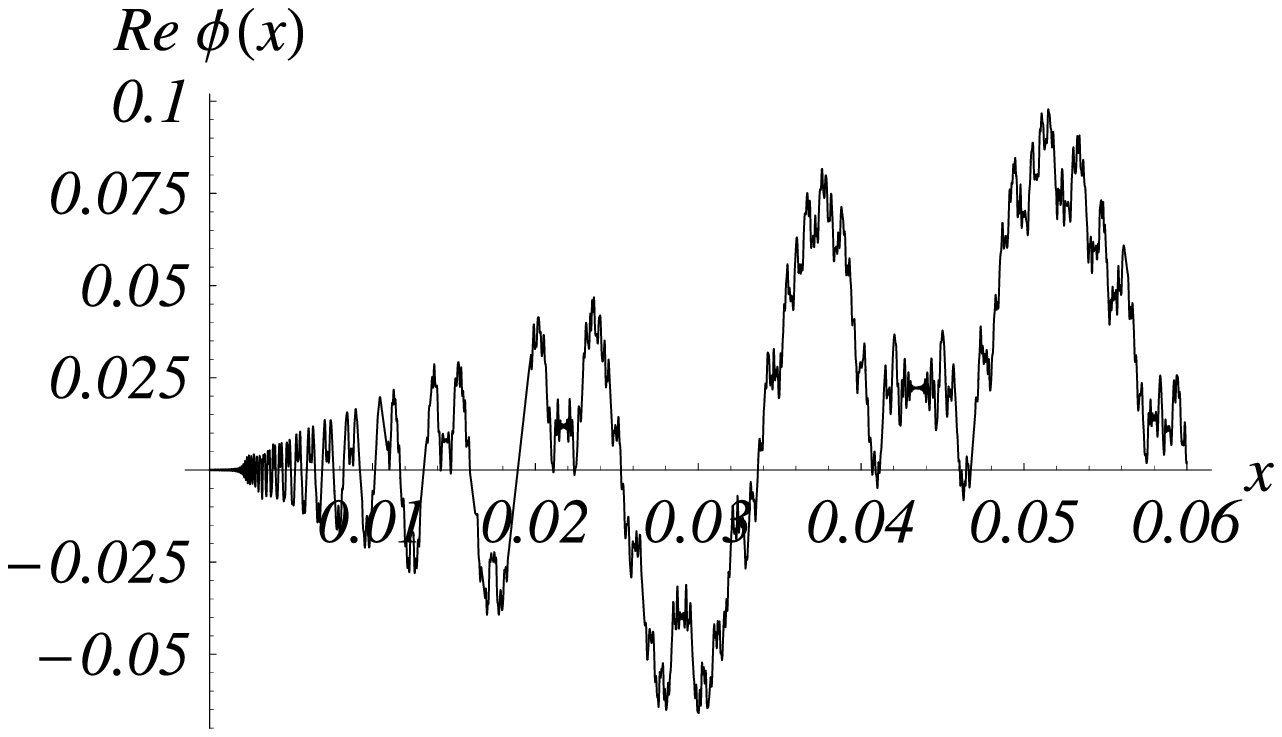}
\end{center}
\caption{Same as figure \ref{weight1big+_3}, but near the origin.} \label{weight1small+_3}
\end{figure}

Modular forms of weight one are another source of continuous, nowhere differentiable functions. The holomorphic function
\begin{equation}\label{weight1}
F(z)\,=\, \frac 12\sum_{(m,n)\in\Z^2} \biggl( e\bigl((m^2+mn+6n^2)z\bigr) - e\bigl((2m^2+mn+3n^2)z\bigr)\biggr)
\end{equation}
is a cuspidal modular form of weight one, automorphic with respect to the subgroup of $SL(2,\Z)$ commonly denoted by $\G_0(23)$. The limit $\,\tau(x)=\lim_{y\to 0^+} F(x+iy)$ exists as a distribution. It has~a continuous, nowhere differentiable first antiderivative, whose real part is graphed in figures \ref{weight1big+_3}-\ref{weight1small+_3}.

Let us describe the content of our paper more closely. In section \ref{autosec} we give quick introduction to automorphic distributions for subgroups of finite index $\G\subset SL(2,\Z)$. We recall the results of \cite{S1} on the regularity behavior of these automorphic distributions in section \ref{regsec}, and then relate these to other known results. Automorphic distributions $\,\tau$ corresponding to Maass forms or modular forms of weight one have continuous first antiderivatives $\phi_\tau$. In section \ref{ratsec} we show that the functions $\Re\phi_\tau$ are non-differentiable at rational points, except in certain special cases. We also explicitly describe the derivatives at rational points when they do exist. Their behavior at irrational points is the subject of section \ref{ratsec}: the $\Re\phi_\tau$ violate certain pointwise H\"older conditions, and that rules out the existence of a derivative. We conclude the paper with a brief discussion of modular forms of weight one half. Our arguments apply almost directly to this case as well, even thought the corresponding automorphic distributions lie in representations of the metaplectic cover of $SL(2,\R)$. Automorphic distributions corresponding to modular forms of weight two have continuous second antiderivatives, which barely miss being differentiable. In the final section, we also present the graph of the imaginary part of the second antiderivative of a certain modular form of weight two.

After this paper was completed, Robert Stanton brought Chamizo's recent paper \cite{Ch} to our attention, which overlaps ours to some extent: Chamizo also produces non-differentiable continuous functions from modular forms, with arguments which are specific to the holomorphic case. According to our view, the phenomena we discuss are caused by automorphy, and can be understood by general arguments.

We are indebted to three colleagues who helped us with certain aspects of this paper: Michael Rubinstein supplied us with the Fourier coefficients of Maass forms, without which we could not have drawn figures 1-3; Henryk Iwaniec and Peter Sarnak enlightened us about the history of the bound (\ref{taufourier3}) for Maass forms.

\section{Automorphic Distributions}\label{autosec}

The automorphic distributions we consider arise as boundary values of mo\-dular forms and Maass forms. We refer the reader to \cite{S1,MS0,MS1} for details about the logical connections between automorphic forms and automorphic distributions. Here we simply recall the technical definitions. Throughout this paper, we use the following notational conventions:
\begin{equation}
\label{GandGamma}
G\ = \ SL(2,\R)\,,\ \ \ \G \subset SL(2,\Z)\ \ \text{a normal subgroup of finite index}.
\end{equation}
The group $G$ and its subgroup $\G$ act on $\R\mathbb P^1=\R\cup\{\infty\}$ by linear fractional transformations.

We write $C^{-\infty}(\R)$ for the space of complex-valued distributions on the real line. According to our convention, distributions are dual to compactly supported smooth measures. Thus functions are special cases of distributions, and distributions ``transform like functions". For $\lambda\in\C$ and $\d\in\Z/2\Z$, we define
\begin{equation}\label{vlambda}
\begin{aligned}
V_{\l,\d}^{-\infty}\ \ = \ \ &\text{vector space of pairs}\ \ (\tau,\tilde\tau)\in C^{-\infty}(\R)\times C^{-\infty}(\R)
\\
&\text{such that}\ \ \tilde\tau(x)=(\sg x)^\d|x|^{\l-1}\tau(-1/x)\ \ \text{for}\,\ x \neq 0\,.
\end{aligned}
\end{equation}
Then $\tau$ determines $\tilde\tau$ except at $x=0$. We shall soon impose a condition that~-- in the cases we are interested in~-- effectively extends $\tilde\tau$ from $\R-\{0\}$ to $\R$. Anticipating this state of affairs, we now tacitly identify each pair $(\tau,\tilde\tau)$ with its first member $\tau$. With that convention, we can describe an action of $G$ on $V_{\l,\d}^{-\infty}$, as follows:
\begin{equation}\label{pilambda}
\begin{aligned}
&\text{for}\ \ g^{-1}\ = \ \begin{pmatrix} a & b \\ c & d   \end{pmatrix}\ \in \ G\,,
\\
&\qquad\bigl(\pi_{\l,\d}(g)\,\tau\bigr)(x)\ = \ (\sg(cx+d))^\d\,|cx+d|^{\lambda-1}\,\tau\left(\frac{ax+b}{cx+d}\right)\,;
\end{aligned}
\end{equation}
strictly speaking, this makes sense only for $cx+d\neq 0$. It can be given meaning even at the missing point by expressing $\tau$ near $x=\infty$ in terms of $\tilde\tau$ near $x=0$. In the special case of $a=d=0$, $b=-c=-1$, $\pi_{\l,\d}(g^{-1})$ simply switches the roles of $\tau$ and $\tilde\tau$. This latter observation implies a formula like (\ref{pilambda}) also for the second member of the pair whose first member is $\pi_{\l,\d}(g^{-1})\tau$. One can check that $\pi_{\l,\d}$ does define a representation of $G$ on $V_{\l,\d}^{-\infty}$, either by a direct computation, or more intelligently, by identifying $V_{\l,\d}^{-\infty}$ with a space of distributions on $G$ -- see \cite{S1}, for example. By definition,
\begin{equation}\label{vgamma}
\bigl(V_{\l,\d}^{-\infty}\bigr)^\G\  = \ \ \text{space of $\G$-invariants in}\ \ V_{\l,\d}^{-\infty}
\end{equation}
is the space of $\G$-{\em automorphic distributions} of type $(\l,\d)$.

We now consider a particular automorphic distribution $\tau\in (V_{\l,\d}^{-\infty})^\G$. Since $\G$ has finite index in $SL(2,\Z)$, there exists a positive integer $N=N(\G)$ such that
\begin{equation}\label{Ngamma}
\begin{pmatrix} 1 & n \\ 0 & 1   \end{pmatrix}\ \in \ \G\ \ \ \Longleftrightarrow\ \ \ n/N \in \Z\,.
\end{equation}
In view of (\ref{pilambda}), the distribution $\tau(x)$ is then periodic of period $N$, and hence has a Fourier expansion
\begin{equation}\label{taufourier}
\tau(x) \ =\ c_0 \ + \ {\sum}_{n\neq 0}\ c_n\, e(nx/N)\qquad \left(\ e(x)\, =_{\text{def}}\, e^{2\pi i x}\  \right)\,.
\end{equation}
To extend $\tau$ across $\infty$, or more precisely, to have $\tau$ determine $\tilde\tau$ completely, we need to give meaning to the distribution
\begin{equation}\label{tildetau}
c_0\,(\sg x)^\d\,|x|^{\lambda-1}\ \ + \,\ \ (\sg x)^\d\,|x|^{\lambda-1}\, {\sum}_{n\neq 0}\ c_n\, e\bigl(-n/(Nx)\bigr)
\end{equation}
even at $x=0\,$. If $\lambda\notin (2\Z+\d)\cap \Z_{\leq 0}\,$,  the distribution $(\sg x)^\d|x|^{\lambda-1}$ can be continued across $x=0$ by analytic continuation in the complex variable $\lambda\,$. From now on we suppose\addtocounter{equation}{1}
\begin{equation*}\label{eitheror1}\tag{\theequation\,a}
c_0\ = \ 0\ \ \ \text{unless}\ \ \ \Re \lambda > 0\,,
\end{equation*}
thus making the first summand in (\ref{tildetau}) well defined at $x=0$. The second summand can be extended across $x=0$ by successive integration by parts, for all values of $\l\in\C\,$; in the terminology of \cite{MS2}, the second summand has a {\em canonical extension} across $0$. When the first summand is extended by virtue of the assumption (\ref{eitheror1}) and the second summand by means of the canonical extension, we say that $\,\tilde\tau$ agrees with its natural extension across $x=0$ or, in terms of $\,\tau$, that $\,\tau$ {\em agrees with its natural extension across $x=\infty$}. We suppose this is the case:
\begin{equation*}\label{eitheror2}\tag{\theequation\,b}
\tau\ \ \text{agrees with its natural extension across}\ \ x=\infty\,.
\end{equation*}
These two conditions imply in particular that the Fourier expansion (\ref{taufourier}) determines $\tau$ not only as distribution on $\R$, but even as element of $(V_{\l,\d}^{-\infty})^\G$.

In order to understand properties of $\,\tau$, we need to work not only with $\,\tau$ itself, but also with its $SL(2,\Z)$-translates. We therefore impose the conditions (\ref{eitheror1},b) not only on $\,\tau$, but on all its translates; we also specifically exclude the case of a constant $\,\tau\in\C-\{0\} \subset (V_{1,0}^{-\infty})^\G$, which we would otherwise have to exclude later, as a trivial counterexample to various statements:
\begin{equation}\label{eitheror}
\begin{aligned}
&\text{the analogues of (\ref{eitheror1}) and (\ref{eitheror2}) hold for all translates}
\\
&\ \ \ \pi_{\l,\d}(\g)\,\tau,\ \g\in SL(2,\Z)\,;\ \text{and $\,\tau\notin\C-\{0\} \subset (V_{1,0}^{-\infty})^\G$}.
\end{aligned}
\end{equation}
We make this a standing assumption throughout our paper. Recall that the group $SL(2,\Z)/\G$ is finite, so (\ref{eitheror}) puts restrictions on only finitely many translates.

\begin{definition}\label{defcuspidal}
An automorphic distribution $\tau\in (V_{\l,\d}^{-\infty})^\G$ satisfying the condition (\ref{eitheror}) is said to be cuspi\-dal at infi\-nity if $c_0=0$; $\,\tau$ is cuspidal (without qualification) if all its $SL(2,\Z)$-translates are cuspidal at infinity.
\end{definition}

Our next statement makes the connection between $\G$-automorphic distributions and more familiar objects. It is not deep. A proof can be extracted from arguments in \cite{S1}. To keep the statement simple, by a {\em holomorphic modular form} for $\G$ we shall mean the datum of two holomorphic modular forms in the usual sense, defined respectively on the upper and the lower half plane. By a {\em Maass form} we mean a Maass form in the strictest sense, i.e., a $\G$-invariant eigenfunction of the $G$-invariant Laplace operator on the upper half plane.

\begin{thm}\label{bijection}
The space of $\,\tau\!\in\! (V_{\l,\d}^{-\infty})^\G\!$ which satisfy the condition (\ref{eitheror}) corresponds bijectively to the space of
\begin{itemize}
\item[{\rm a)}]
cuspidal Maass forms for $\G$ with eigenvalue $\frac 14(1-\l^2)\geq \frac 14$, in case $\l\in i\R$ and $\,\d=0\,$;
\item[{\rm b)}]
cuspidal Maass forms for $\G$ with eigenvalue $\frac 14(1-\l^2)< \frac 14$, in case $\,-1<\l<0\,$ and $\,\d=0\,$;
\item[{\rm c)}]
square-integrable Maass forms for $\G$ with eigenvalue $\frac 14(1-\l^2)< \frac 14$, in case $\,0<\l<1\,$ and $\,\d=0\,$;
\item[{\rm d)}]
cuspidal odd-weight Maass forms for $\G$, of any given odd weight, with eigenvalue $\frac 14(1-\l^2)>\frac 14$, in case $\l\in i(\R-\{0\})$ and $\,\d=1\,$;
\item[{\rm e)}]
cuspidal holomorphic modular forms of weight $k\geq 1$, in case $\l = 1-k$ and $\,\d\equiv k \pmod{2}$.
\end{itemize}
In all other cases the space of all such $\,\tau$ reduces to $0\,$.
\end{thm}

Some or all of the Maass forms in c) may be cuspidal, depending on whether or not the corresponding $\,\tau$ are cuspidal; in all other cases, cuspidality follows from the standing hypothesis (\ref{eitheror}). In the situations a) and d), both $\l$ and $-\l$ determine the same spaces of Maass forms. The standard intertwining operator $V_{-\l,\d}^{-\infty}\simeq V_{\l,\d}^{-\infty}$ explains this coincidence on the level of automorphic distributions. The intertwining operator also relates the automorphic distributions in b) to the cuspidal cases in c).

In effect, our condition (\ref{eitheror}) rules out Eisenstein series, as well as the images of the non-cuspidal automorphic distributions $\,\tau$ in c) under the standard intertwining operator. It would not be difficult to extend our discussion also to these two cases; they are less interesting, and would make various statements more involved. We should point out that Selberg's eigenvalue conjecture predicts the non-existence of non-zero Maass forms in the situations b) and c) if $\G$ is a congruence subgroup. For certain non-congruence subgroups, Maass forms of this type~-- even cuspidal Maass forms~-- are known to exist \cite{Se}.

\section{Regularity properties}\label{regsec}

Recall that a function $f\in C(\R)$ is said to be H\"older continuous of index $\alpha$, $0<\alpha\leq 1$, if
\begin{equation}\label{hoelder}
|f(x)-f(y)|\ < \ C|x-y|^\alpha \ \ \text{for all}\ \ x,\, y\in \R\,,
\end{equation}
for some constant $C>0$ which can be chosen locally uniformly in $x,\,y$. The functions that come up in this paper are generally periodic, in which case $C$ can be chosen independently of $x$ and $y$. As in \cite{S1}, we define spaces the spaces $C^\alpha(\R)\subset C^{-\infty}(\R)$, $\alpha \in \R$, as follows: for $0<\alpha<1$, $C^\alpha(\R)$ is the space of $\alpha$-H\"older continuous functions; $C^0(\R)$ is the space of continuous functions; and we extend the definition to functions and distributions so that
\begin{equation}\label{calpha}
C^\alpha(\R)\ = \ \textstyle\frac{d\ }{dx}\,C^{\alpha+1}(\R)\ \ \ \text{for all $\,\alpha\in\R$}\,.
\end{equation}
This results in the usual definition of the space $C^k(\R)$ when $k=\alpha\in\Z_{\geq 0}$; for $k\in\Z_{>0}$, $C^{-k}(\R)$ is the space of distributions expressible as $k$-th derivative of a continuous function. Further notation:
\begin{equation}\label{c<alpha}
C^{<\alpha}(\R)\ = \ {\cap}_{\beta<\alpha}\ C^\beta(\R)
\,,\qquad C^{>\alpha}(\R)\ = \ {\cup}_{\beta>\alpha}\ C^\beta(\R) \,.
\end{equation}
Thus $C^\gamma(\R)\subset C^{>\beta}(\R)\subset C^{\beta}(\R)\subset C^{<\beta}(\R)\subset C^\alpha(\R)$ whenever $\alpha<\beta<\gamma$.

\begin{thm}[\cite{S1}]\label{regthm}
Under the standing hypothesis (\ref{eitheror}) on $\tau\in (V_{\l,\d}^{-\infty})^\G$,
\begin{itemize}
\item[{\rm a)\ }]
$\tau\in C^{\l - 1}(\R)\,$ if $\,\tau$ is non-cuspidal;
\item[{\rm b)\ }]
$\tau\in C^{\frac {\Re \l-1}2}(\R)\,$ if $\,\tau$ is cuspidal and $\,\lambda-1\notin 2\Z\,$;
\item[{\rm c)\ }]
$\tau\in C^{<\frac {\l-1}2}(\R)\,$ if $\,\tau$ is cuspidal and $\,\lambda-1\in 2\Z\,$.
\end{itemize}
\end{thm}
\medskip

In particular, $\,\tau$ has a continuous anti-derivative if it corresponds to a square integrable Maass form with $\d=0$ and $0<\l<1$, to any cuspidal Maass form, possibly of odd weight, or to a cuspidal holomorphic modular form of weight one.

The next statement asserts the equivalence of various regularity
criteria for a periodic distribution $\,\tau$, as in
(\ref{taufourier}). The origin of $\,\tau$ will not matter now~--
in other words, $\,\tau$ need not be $\G$-automorphic. Regularity
is not affected by adding a constant, nor by scaling of the
variable. We may therefore suppose that $c_0=0$ and $N=1$:
\begin{equation}\label{taufourier1}
\tau(x)\ \ = \ \ {\sum}_{n\neq 0}\ c_n\, e(nx)\,.
\end{equation}
Our statement also involves
\begin{equation}\label{taufourier2}
\tau_+(x)\ \ = \ \ {\sum}_{n> 0}\ c_n\, e(nx)\,,\ \ \ \tau_-(x)\ \ = \ \ {\sum}_{n< 0}\ c_n\, e(nx)\,.
\end{equation}

\begin{thm}\label{regthm2}
{\rm 1)}\ \ For $k\in\N$ and $\alpha\in\R$, with $\a\leq k$, the following conditions are equivalent:
\begin{itemize}
\item[{\rm a)\ }]
$\tau\in C^{<\alpha}(\R)\,$;
\item[{\rm b)\ }]
$\tau_+\in C^{<\alpha}(\R)\,$ and $\,\tau_-\in C^{<\alpha}(\R)$;
\item[{\rm c)\ }]
for every $\e>0$, there exists a constant $C=C(\e)$ such that\\ $\bigl|{\sum}_{|n|\leq N}\,c_n \,n^k\, e(nx)\bigr|\leq C\,N^{\e+k-\a}$ as $N\to\infty$, uniformly in $x$;
\item[{\rm d)\ }]
for every $\e>0$, there exists a constant $C=C(\e)$ such that both\\ $\bigl|{\sum}_{n=1}^N\, c_n\,n^k e(nx)\bigr|\leq C N^{\e+k-\a}$ and $\,\bigl|{\sum}_{n=1}^N \, c_{-n}\, n^k e(-nx)\bigr|\leq C N^{\e+k-\a}$ as $N\to\infty$, uniformly in $x$.
\end{itemize}
{\rm 2)}\ \ For $\a\in\R\,$ and $s\in\C\,$, the condition {\rm b)} is equivalent to
\begin{itemize}
\item[{\rm e)\, }]
${\sum}_{n>0}\, \frac{c_n}{n^s}\, e(nx) \in C^{<(\a+\Re s)}(\R)$ and $\,{\sum}_{n>0}\, \frac{c_{-n}}{n^s}\,e(-nx)  \in C^{<(\a+\Re s)}(\R)$.
\end{itemize}
\end{thm}

The proof uses only standard tools, but we have not been able to find these statements assembled in the literature. Before turning to the proof we discuss some of the implications.

\begin{cor}[S. Bernstein's Theorem \cite{Be}]
If $\,\tau\in C^{>\alpha}(\R)\,$,
\[
{\sum}_{n\neq 0}\ |n|^{2\a}\,|c_n|^2\ < \ \infty\ \ \ \text{and}\ \ \  {\,\sum}_{n\neq 0}\ |n|^{\a-1/2}\,|c_n|\ < \ \infty\,.
\]
\end{cor}

\begin{proof}
According to the theorem, ${\sum}_{n\neq 0}\,c_n\, n^\a\,e(nx) \in C^{>0}(\R/\Z)$, hence
\begin{equation}\label{bernstein1}
{\sum}_{n\neq 0}\,|n|^{2\a}\,|c_n|^2\ = \ \int_{\R/\Z}\ \bigl|{\sum}_{n\neq 0}\,c_n\, n^\a\,e(nx)\bigr|^2\,dx \ < \ \infty\,.
\end{equation}
For the second assertion, we appeal to the Cauchy-Schwartz inequality:
\begin{equation}\label{bernstein2}
\begin{aligned}
{\sum}_{n\neq 0}\,|n|^{\a - 1/2}\,|c_n|\ &= \ {\sum}_{n\neq 0}\,\bigl(|n|^{\a +\e}\,|c_n|\bigr)\,|n|^{-1/2-\e}
\\
&\leq \ \bigl({\sum}_{n\neq 0}\,|n|^{2\a+2\e}\,|c_n|^2\bigr)^{\frac 12}\bigl( {\sum}_{n\neq 0}\,|n|^{-1-2\e} \bigr)^{\frac 12}\,;
\end{aligned}
\end{equation}
the right hand side is finite for every sufficiently small $\e>0$, because the hypothesis permits us to enlarge the constant $\a$ in (\ref{bernstein1}) slightly.
\end{proof}

We combine this classical theorem of S. Bernstein with theorem \ref{regthm}, to obtain a regularity condition of Sobolev type: if the automorphic distribution $\,\tau\in (V_{\l,\d}^{-\infty})^\G$ satisfies the condition (\ref{eitheror}),
\begin{equation}\label{bernrez}
{\sum}_{n}\,|n|^{\a}\,|c_n|^2\ < \ \infty\,,\ \ \text{provided}\ \ \begin{cases} \a<\Re\l-1\ &\text{if $\,\tau$ is cuspidal} \\ \a< 2(\l-1) &\text{otherwise;}\end{cases}
\end{equation}
recall that our hypotheses imply $0<\l<1$ when $\,\tau$ fails to be cuspidal. This, in effect, is Bernstein-Reznikoff's regularity theorem \cite{BR}. Note, however, that the Sobolev bounds do not imply the H\"older bounds of theorem \ref{regthm}.

One can define a fractional derivative of order $\b>0$ for any periodic distribution $\,\tau$ without constant term, as in (\ref{taufourier1})\,:
\begin{equation}\label{fracderiv}
\tau^{(\b)}(x)\, = \, (2\pi)^\b\left(e^{\frac{\b \pi i}{2}}\,{\sum}_{n>0}\, c_n n^\b e(nx) +  e^{\frac{3\b \pi i}{2}}\,{\sum}_{n<0}\, c_n|n|^\b e(nx)\right).
\end{equation}
As one consequence of the theorem, $\,\tau\in C^{<\a}(\R)$ implies $\,\tau^{(\b)}\in C^{<\a-\b}(\R)$, and conversely. A slightly stronger result is due to Hardy and Littlewood \cite{Zy}: $\,\tau\in C^{\a}(\R)\,$ if and only if $\,\tau^{(\b)}\in C^{\a-\b}(\R)$, unless $\a\in\Z$ or $\a-\b\in\Z$.

The traditional approach to Maass forms involves a Fourier expansion on the upper half plane, in terms of Bessel functions. The resulting Fourier coefficients $a_n$ are related to the Fourier coefficients $c_n$ in (\ref{taufourier}),
\begin{equation}\label{taufourier2.5}
a_n\ \ = \ \ |n|^{\l/2}\,c_n\,.
\end{equation}
The coefficients $a_n$ are used to define the $L$-function of a Maass form, which results in a functional equation with symmetry $s \mapsto 1-s$. In the context of modular forms of weight $k$, one usually defines the $L$-function in terms of the $c_n$, with symmetry $s\mapsto k-s$. Sometimes the $a_n$ are used instead, especially if one wants to discuss Maass forms and modular forms on an equal footing; this is the so-called ``Langlands normalization". Theorems \ref{regthm} and \ref{regthm2} imply a bound on the $a_n$\,: for every $\e>0$,
\begin{equation}\label{taufourier3}
{\sum}_{|n|\leq N}\ a_n\,e(nx)\ = \ O(N^{1/2+\e})\,,\ \ \text{provided $\,\tau$ is cuspidal}\,.
\end{equation}
Conversely this estimate implies $\,\tau\in
C^{<\frac{\Re\l-1}{2}}(\R)$, which is within a hair's breadth of
the H\"older estimate of theorem \ref{regthm}. The estimate
(\ref{taufourier3}) is completely standard in the case of modular
forms; for Maass forms, this estimate appears in several places,
but we are aware of only one completely satisfactory proof
\cite{Har} -- besides the argument we have just given, of course.
The same estimate is expected to hold for the coefficients $a_n$
of $L$-functions of cuspidal automorphic representations for
$GL(n)$ -- this is the so-called cancellation conjecture. For
$GL(3)$, an estimate halfway between the cancellation conjecture
and the obvious bound $O(N)$ is the proven in
\cite{cancellation2}.\bigskip

\noindent{\it Proof\/}\ of theorem \ref{regthm2}. We shall use the notational convention $f_N(x) = O(N^\a)$ for a family of functions $f_N$ to signify an estimate $|f_N(x)|\leq CN^\alpha$ that holds uniformly in $x$, for all but finitely many $N$. The notation $f_N(x) = O_\e(N^{\a+\e})$ signifies the same kind of bound, for every $\e>0$, with $C=C_\e$ allowed to depend on $\e$. Recall the partial summation formula:
\begin{equation}\label{partialsum}
{\sum}_{n=1}^N\, a_n\,b_n \ = \ {\sum}_{n=1}^{N-1} \bigl(b_n - b_{n+1} \bigr)\biggl( {\sum}_{k=1}^n a_k\biggr)\ + \ b_N\, {\sum}_{n=1}^N\, a_n\,.
\end{equation}
We apply this with $a_n=c_ne(nx)$ and $b_n= n^{-s}$, for some fixed $s\in\C$. Since $n^{-s} - (n+1)^{-s} = O(n^{-\Re s-1})$,
\begin{equation}\label{partialsum1}
{\sum}_{n=1}^N\, c_n\, e(nx)\ = \ O(N^\a)\ \ \Longleftrightarrow\ \ {\sum}_{n=1}^N\, \frac{c_n}{n^s}\, e(nx) \ = \ O(N^{\a-\Re s})\,,
\end{equation}
provided $\,\a> 0$ and $\,\a-\Re s > 0$\,; one direction follows directly from (\ref{partialsum}), and the other comes ``for free", since $\,\a$ and $\,\a-\Re s$ play symmetric roles.

For the equivalence between c) and d), and also to relate these conditions to a) and b), we need the Dirichlet kernel and the ``half Dirichlet kernel"
\begin{equation}\label{dirichlet}
\begin{aligned}
D_N(x)\ &= \ {\sum}_{|n|\leq N}\, e(nx) \, = \, \frac{\sin\bigl((2N+1)\pi x\bigr)}{\sin(\pi x)}\,,
\\
P_N(x)\, &= \, {\sum}_{0\leq n\leq N}\, e(nx) \, = \, \frac{e\bigl((N+1)x\bigr)-1}{e(x)-1}\, = \, e(\textstyle\frac {Nx}2)\,D_{N/2}(x)\,;
\end{aligned}
\end{equation}
one can make sense of the expression for $D_N(x)$ as a quotient of sine functions even when $N$ is a half-integer, e.g., for $|x|<1$, thus giving meaning to the relation between the two kernels for all $N\in\N$. We shall need to know that
\begin{equation}\label{dirichlet1}
\|D_N(x)\|_1\ = \ O(\log N)\ \ \text{and}\ \ \|P_N(x)\|_1\ = \ O(\log N)\,;
\end{equation}
this can be seen by integrating $|D_N(x)|$ first over the interval $[0,1/N]$ -- on which $|D_N(x)|$ is bounded by a multiple of $N$~-- then over the interval $[1/N,1/2]$~-- on which $\sin(\pi x)$ is bounded from below by a multiple of $x$~-- and finally using the symmetry about $x=1/2$. The estimate for $D_N$ then implies the same type of estimate for $P_N$.

Obviously d) implies c). For the converse, note that convolution with $P_N(x)$ maps the partial sum $\sum_{|n|\leq N}c_n e(nx)$ to the sum $\sum_{0<n\leq N}c_n e(nx)$. This process increases the supremum norm at most by a factor equal to the $L^1$ norm of $P_N$, i.e., at most by the factor $\log N$. This shows:
\begin{equation}\label{dirichlet2}
{\sum}_{n=-N}^N\, c_n\, e(nx)\ = \ O(N^\a)\ \ \Longrightarrow\ \ {\sum}_{n=1}^N\, c_n\, e(nx) \ = \ O(N^\a \log N)\,.
\end{equation}
Since $\log N = O(N^\e)$ for any $\e >0$, we may conclude that c) and d) are equivalent.

Because of (\ref{partialsum1}), the assertion d) for $\alpha$ and
$k$ implies the same assertion for the same $\a$ and any other
choice of $k$, as long as $\a$ remains less than or equal to $k$.
In particular, the criterion d) is compatible with
differentiation: if $\,\tau_+$ and $\,\tau_-$ satisfy d), the
derivatives $\,\tau_+'$, $\,\tau_-'$ satisfy the same condition
with $\a-1$ and $k-1$ in place of $\a$ and $k$, and conversely.
The conditions a) and b) are similarly compatible with
differentiation; cf. (\ref{calpha}). Hence, from now on, we may
suppose that $0< \a \leq 1$ and $k=1$. We already know the
equivalence of c) and d). Since b) trivially implies a), we only
need to deduce c) from a) and b) from d); once a)\,--\,d) are
known to be equivalent, 2) follows formally from
(\ref{partialsum1}).

With $0< \a \leq 1$ and $k=1$ as we are assuming, we express the relevant partial sum as an integral against the derivative of the Dirichlet kernel,
\begin{equation}\label{dirichlet3}
\begin{aligned}
2\pi i\ {\sum}_{|n|\leq N}\,c_n\,n\,e(nx)\ \ &= \ \ \int_{\R/\Z} \tau(x-t)\,D'_N(t)\,dt
\\
&= \ \ \int_{\R/\Z} \bigl( \tau(x-t)\, - \tau(x)\bigr) \,D'_N(t)\,dt\,;
\end{aligned}
\end{equation}
at the second step we have used the fact that $D'_N(x)$ has a Fourier series without constant term. We split the integration into three parts, from $-1/N$ to $1/N$, from $-1/2$ to $-1/N$, and finally from $1/N$ to $1/2$. On the first, $|D'_N(x)|$ can be bounded by a multiple of $N|x|^{-1}$, on the other two by a multiple of $x^{-2}$, independently of of $N$. We choose $\e>0$ so that $0<\a-\e<1$. Since $|\tau(x-t)\, - \tau(x)|\leq C\,|t|^{\alpha-\e}$ uniformly in $x$ and $t$, (\ref{dirichlet3}) implies
\begin{equation}\label{dirichlet4}
\begin{aligned}
\bigl|\, {\sum}_{|n|\leq N}\,c_n\,n\,e(nx)\,\bigr|\ \leq \ C_1 \int_{-1/N}^{1/N} &N\,|t|^{\a-\e-1}\,dt\ +
\\
&+ \ 2\,C_2 \int_{1/N}^{1/2} |t|^{\a-\e-2}\,dt\,,
\end{aligned}
\end{equation}
with $C_1$, $C_2$ depending on $\,\tau$, $\a$ and $\e$, but not on $N$. Since $\e$ is arbitrary except for the restriction $0<\e<\a$, the left hand side has order of growth $O_\e(N^{\e+1-\a})$. That establishes the implication a) $\Longrightarrow$ c).

Finally we start with the hypothesis d) with $0<\a\leq 1$, $k=1$, and $\e>0$ chosen so that $\a-\e>0$. We apply (\ref{partialsum}) with $a_n=c_n n e(nx)$, $b_n=n^{-1}$, summing first from $1$ to $N$, then from $1$ to $M<N$, then taking the difference:
\begin{equation}\label{dirichlet5}
\begin{aligned}
{\sum}_{n=M}^N\,c_n\,e(nx)\ = \ {\sum}_{n=M}^{N-1}\biggl( {\sum}_{\ell =1}^n\,\ell\, c_\ell\,e(\ell x) \biggr)\bigl(n^{-1} - (n+1)^{-1}\bigr)
\\
+ \ N^{-1}\,{\sum}_{\ell =1}^N\, \ell\,c_\ell \,e(\ell x)\ - \ M^{-1}\,{\sum}_{\ell =1}^M \,\ell\, c_\ell \,e(\ell x)\,.
\end{aligned}
\end{equation}
Since ${\sum}_{\ell =1}^n \, \ell\,c_\ell\,e(\ell x) = O_\e(n^{\e + 1 - \a})$, $\,\frac 1n -\frac 1{n+1}=O(n^{-2})$, and $\,M<N$,
\begin{equation}\label{dirichlet6}
{\sum}_{n=M}^N\,c_n\,e(nx)\ = \ O_\e(M^{\e-\a})\,,
\end{equation}
independently of $N$. In particular the series for $\,\tau_+$ converges uniformly, and $\,\tau_+$ is continuous. In proving the H\"older bound (\ref{hoelder}), we may as well suppose $|x-y|<1$. We choose the integer $M$ so that $M<|x-y|^{-1}\leq M+1$. Then
\begin{equation}\label{dirichlet7}
\begin{aligned}
\tau_+(x) - \tau_-(y)\ = \ &{\sum}_{n=1}^{M-1}\,c_n\,e(nx)\ - \ {\sum}_{n=1}^{M-1}\,c_n\,e(ny) \ +
\\
&+ \ {\sum}_{n=M}^\infty\,c_n\,e(nx)\ - \ {\sum}_{n=M}^\infty\,c_n\,e(ny)\,.
\end{aligned}
\end{equation}
Letting $N$ tend to $\infty$ in (\ref{dirichlet6}), one sees that the last two sums each have order $O_\e(M^{\e-\a})$. The difference of the first two terms on the right hand side can be bounded by the product of $|x-y|$ times the supremum of
\begin{equation}\label{dirichlet8}
\frac {d\ }{dx}\,{\sum}_{n=1}^{M-1}\,c_n\,e(nx)\ = \ 2\pi i\, {\sum}_{n=1}^{M-1}\,n\,c_n\,e(nx) \ = \ O_\e(M^{1+\e-\a})\,.
\end{equation}
But $|x-y| \sim M^{-1}$, letting us conclude $|\tau_+(x)-\tau_+(y)|= O_\e(|x-y|^{\a-\e})$ as $x\to y$. In other words, $\,\tau_+\in C^{<\a}(\R/\Z)$. Similarly $\,\tau_-\in C^{<\a}(\R/\Z)$, since $\,\tau_+$ and $\,\tau_-$ play symmetric roles. Thus d) $\Longrightarrow$ b).
\bx

\section{Behavior near rational points}\label{ratsec}

In this and the next section, $\tau\in (V_{\l,\d}^{-\infty})^\G$ denotes an automorphic distribution, subject to the usual hypothesis (\ref{eitheror}), which can be expressed as the first derivative of a continuous function. According to theorem \ref{regthm}, that is the case when the parameter $(\l,\d)\,$ satisfies either of the following conditions:
\begin{equation}\label{contantider}
\begin{aligned}
{\rm a)\ \ } &\d =0 \ \ \text{or}\ \ \d = 1\,,\ \ \ \l \in i\,\R\,,
\\
{\rm b)\ \ }&\d =0\,,\ \ \ -1 < \l < 0\ \ \text{or}\ \ 0 < \l < 1\,.\qquad\qquad\qquad\qquad\qquad\qquad
\end{aligned}
\end{equation}
This covers all cases in theorem \ref{bijection} except holomorphic modular forms of weight $k\geq 2$. In fact, theorem \ref{regthm} is sharp; when $\,\tau$ corresponds to a holomorphic modular form of weight at least $2$, it does not have a continuous anti-derivative. We therefore suppose that one of the two conditions (\ref{contantider}\,a,b) holds.

A periodic distribution without constant term has a distinguished anti-derivative, periodic of the same period, and also without constant term. In the situation we are considering, that means
\begin{equation}\label{contantider1}
\tau(x)\ = \ c_{\tau,0}\ + \ \phi'_{\tau}(x)\,,\ \ \ \phi_\tau \in C^0(\R/N\Z)\,,\ \ \ \int_0^N \phi_\tau(x)\,dx\ = \ 0\,;
\end{equation}
the constant $c_{\tau,0}$ is the constant term of the Fourier series (\ref{taufourier}). We incorporate $\,\tau$ into the notation as a subscript because we also need the analogous expressions $\,\pi_{\l,\d}(\g)\,\tau(x)=c_{\g\tau,0}+\phi'_{\g\tau}(x)$ for the various translates $\,\pi_{\l,\d}(\g)\,\tau$, $\,\g\in SL(2,\Z)$. For $k\geq 1$, we let $\phi^{(-k)}_{\g\tau}\in C^k(\R/N\R)$ denote the $k$-th antiderivative of $\phi^{(0)}_{\g\tau}=\phi_{\g\tau}$, normalized by the requirement that its Fourier series have zero constant term.

For the next statement, we fix a rational number $p/q$, expressed as the quotient of relatively prime integers, with $q>0\,$; for emphasis,
\begin{equation}\label{pqcond}
p,\,q\in\Z\,,\ \ (p,q) = 1\,,\ \ q>0\,.
\end{equation}
We can then choose $r,s\in\Z\,$ so that $\,pr-qs=1$, and define\begin{equation}\label{gammapq}
\g\ \ = \ \ \begin{pmatrix} r & -s \\ -q & p \end{pmatrix}\ \in \ SL(2,\Z)\,.
\end{equation}
Note that $\,\g$ maps the point $p/q$ to $\infty$.

\begin{lem}\label{asympexplemma}
Under the hypotheses just stated, for $x\in \R$ and $n\geq 0$,
\[\begin{aligned}
\phi_\tau(x&) \,-\,  \phi_\tau(p/q) \ \ =\ \  \frac{c_{\tau,0}}{q}\,(p-qx)\,\ - \ \frac{c_{\g\tau,0}}{\l\,q} \bigl(\sg(p-qx)\bigr)^{\d+1}|p-qx|^\l
\\
&+\, {\sum}_{k=0}^n\, q^k\bigl(\sg(p-qx)\bigr)^{\d+k}\left( {\prod}_{1\leq j\leq k} (\l+j) \right)\phi_{\g\tau}^{(-k)}(\g\,x)|p-qx|^{\l+k+1}
\\
&-\ q^{n+1}\bigl(\sg(p-qx)\bigr)^{\d+n}\left( {\prod}_{0\leq j\leq n} \,(\l+j+1) \right)\ \times
\\
&\qquad\ \ \times \ \int_{\sg(p-qx)\g x}^{+\infty} \bigl(qt+r\sg(p-qx)\bigr)^{-\l -n-2} \,\phi_{\g\tau}^{(-n)}(\sg(p-qx)t)\,dt\,.
\end{aligned}
\]
\end{lem}

This formula appears already implicitly in \cite{S1}, but only special cases are written out there in detail. If $\l=0$, (\ref{eitheror}) implies the vanishing of $c_{\g\tau,0}$; in that case, the product $\l^{-1}c_{\g\tau,0}$ is to be interpreted as $0$, of course.

\begin{proof}
The expression $(qt+r\sg(p-qx))$ is strictly positive on the interval of integration -- recall that $q>0$ by assumption. Since $\re\l>-1$, the integral converges comfortably, in fact
\begin{equation}\label{intbound}
\begin{aligned}
&\left|\int_{\sg(p-qx)\g x}^{+\infty} \bigl(qt+r\sg(p-qx)\bigr)^{-\l -n-2} \,\phi_{\g\tau}^{(-n)}(\sg(p-qx)t)\,dt\,\right|\ \leq
\\
&\qquad\qquad\leq\ q^{-1}\,(\Re\l+n+1)^{-1}\,|p-qx|^{\Re\l+n+1}\, \operatorname{max}|\phi_{\g\tau}^{(-n)}(x)|\,.
\end{aligned}
\end{equation}
Both sides of the identity vanish for $x=p/q$; here again we are using the fact that $\Re\l\leq 0$ implies $c_{\g\tau,0}=0$. It therefore suffices to equate the derivatives. Differentiation of the identity results in the equation $\,\tau(x)= (\sg(p-qx))^\d |p-qx|^{\l-1}(\pi_{\l,\d}(\g))\tau(\g x)$, or equivalently in the tautological equation $\,\tau = \,\tau$.
\end{proof}

The estimate (\ref{intbound}) only depends the boundedness of the function $\phi_{\g\tau}^{(-n)}$. Using the periodicity and vanishing of the constant term, one can strengthen the estimate by one additional power of $|p-qx|$, to $O(|p-qx|^{\Re\l+n+2})$. Hence
\begin{equation}\label{asymptoticexpension}
\begin{gathered}
\text{the identity in lemma (\ref{asympexplemma}) constitutes an asymptotic}
\\
\text{expansion for $\,\phi_\tau(x)\,$ as $\,x \to p/q$\,}.
\end{gathered}
\end{equation}
The strengthened estimate depends on the integer $N$, among other
things. One can establish (\ref{asymptoticexpension}) also more
easily by noting that the right hand side of the estimate
(\ref{intbound}) also bounds the last term in the sum in lemma
\ref{asympexplemma}, up to a factor involving $\l$ and $q^n$. In
particular, taking $n=0$, we find

\begin{cor}\label{n=1} There exists a constant $C>0$, depending on the maximum absolute value of the (finitely many) functions $\,\phi_{\g_j\tau}^{(-1)}$, $\g_j\in SL(2,\Z)$, and on $\l$, but not otherwise on $\tau$, $\G$, or $p$ and $q$, such that
\[\begin{aligned}
\biggl| \, \phi_\tau(x) \,-\,  \phi_\tau&(p/q)\, - \, \frac{c_{\tau,0}}{q} (p-qx)\, + \, \frac{c_{\g\tau,0}}{\l\,q} \bigl(\sg(p-qx)\bigr)^{\d+1}|p-qx|^\l
\\
& -\,\bigl(\sg(p-qx)\bigr)^{\d} |p-qx|^{\l+1}\phi_{\g\tau}(\g\,x)\, \biggr|\ \leq \ C\,q\,|p-qx|^{\Re\l+2}\,.
\end{aligned}
\]
\end{cor}

The $\,\phi_{\g_j\tau}^{(-1)}$ are periodic of period $N$, with zero constant term. That makes it possible to bound these functions in terms of the $\,\phi_{\g_j\tau}$. In other words, the constant $C$ can be made to depend on the maximum absolute value of the $\,\phi_{\g_j\tau}$ instead, but it then also depends on $N$.

Since $\Re\l+2>1$ by assumption, the differentiability of $\phi_\tau$ at $x=p/q$ is governed completely by the last two terms between the absolute value bars in (\ref{n=1}). The first of the two prevents differentiability unless $c_{\g\tau,0}=0$. In analyzing the second term, recall that $\phi_{\g\tau}$ is periodic of period $N$, with zero constant term. Since
\begin{equation}\label{closetopoverq}
\g x\ = \ \frac{rx-s}{p-qx} \ \sim \ -\,q^{-2}\,\left( x - \frac pq\right)^{-1}\ \ \ \text{as}\ \ x \to \frac pq\,,
\end{equation}
$\phi_{\g\tau}(\g x)$ assumes all values of $\phi_{\g\tau}$ infinitely often near $p/q$, with an approximate spacing of $N(p-qx)^2$ between successive cycles. When $\l$ has a non-zero imaginary part, the phase of $|p-qx|^{\l+1}$ goes through a complete cycle over intervals of approximate length $2\pi|x-p/q|(\Im\l)^{-1}$. Since $(p-qx)^2 \ll |x-p/q|$ for $x\to p/q$,
\begin{equation}\label{closetopoverq1}
\begin{aligned}
&\lim\operatorname{sup}_{x\to p/q}\left(\frac{\Re\left(|p-qx|^{\l+1}\phi_{\g\tau}(x)\right)}{\left|p-qx   \right|^{\Re\l+1}}\right)\ = \ \operatorname{max}|\phi_{\g\tau}(x)|
\\
&\qquad\qquad= \ -\,\lim\operatorname{inf}_{x\to p/q}\left(\frac{\Re\left(|p-qx|^{\l+1}\phi_{\g\tau}(x)\right)}{\left|p-qx   \right|^{\Re\l+1}}\right),
\end{aligned}
\end{equation}
and the analogous identities hold for the imaginary part. An even simpler argument gives a similar conclusion when $\l$ is real. In any case, neither the real part nor the imaginary part of the term under consideration can be differentiable unless $\l>0$ -- or unless $\phi_{\g\tau}(\g x)\equiv 0$, in which case $\,\tau\equiv 0$, too.

\begin{cor}\label{diffatrationals} Suppose $\,\tau$ does not vanish identically. Then $\Re\phi_\tau(x)$ and $\Im\phi_\tau(x)$ fail to be differentiable at $x=p/q$, except in the following situation: $\l>0$ and $\Re c_{\g\tau,0}=0$, respectively $\Im c_{\g\tau,0}=0$. When that is the case, $\Re\phi'_\tau(p/q)=-\Re c_{\tau,0}$, respectively $\Im\phi'_\tau(p/q)=-\Im c_{\tau,0}$. In particular, if $\l>0$ and if $\,\tau$ is cuspidal, $\phi_\tau(x)$ is differentiable at all rational points, with derivative $0$.
\end{cor}

The notation $c_{\g\tau,0}$ might suggest that this quantity depends on $\,\g$ -- in fact, it depends only on $p$ and $q$, as can be checked by going back to the definition of $c_{\g\tau,0}$ as the constant term in the Fourier series for $\pi_{\l,\d}(\g)\tau$.

If $\Re\l=0$, our argument shows that $\phi_\tau(x)$ satisfies a pointwise H\"older condition at $x=p/q$, with H\"older index $1$, i.e.
\begin{equation}\label{closetopoverq2}
 \phi_\tau(x)-\phi_\tau(p/q)   \ \ = \ \ O(\,|x-p/q|\,) \ \ \ \text{as}\ \ x \to p/q\,.
\end{equation}
That by itself would not rule out differentiability at $p/q$. In the other cases of non-differentiability, we have established slightly more than stated in the corollary: if $\l>0$ and $\Re c_{\g\tau,0}\neq 0$, (\ref{n=1}) implies
\begin{equation}\label{closetopoverq3}
{\lim}_{x\to p/q}\ |x-p/q|^{-\alpha}\left|\Re\bigl( \phi_\tau(x)-\phi_\tau(p/q)  \bigr)   \right| \ = \ \infty\ \ \text{if}\ \ \alpha>\lambda\,.
\end{equation}
Since $\l<1$, this does rule out differentiability in a quantitative manner. On the other hand, if $\l<0$, (\ref{closetopoverq1}) implies
\begin{equation}\label{closetopoverq4}
{\lim}_{x\to p/q}\ |x-p/q|^{-\alpha}\left|\Re\bigl( \phi_\tau(x)-\phi_\tau(p/q)  \bigr)   \right| \ = \ \infty\ \ \text{if}\ \ \alpha>1+\lambda\,,
\end{equation}
again ruling out differentiability.

\section{Behavior near irrational points}\label{irratsec}

We continue with the hypotheses of the previous section, in particular (\ref{contantider}), which ensures the existence of a continuous anti-derivative. Any irrational $x_0\in\R$ can be approximated arbitrarily closely by a sequence of rational numbers $p/q$ such that
\begin{equation}\label{generic}
|\,p\,-\, qx_0\,|\ < \ q^{-1}\,.
\end{equation}
Generically the exponent $\,-1$ cannot be improved, in which case one says that $x_0$ has {\em irrationality measure $2$}. There do exist irrational numbers $x_0$ with an approximating sequence of rationals $p/q$ such that
\begin{equation}\label{nongeneric}
|\,p\,-\, qx_0\,|\ < \ q^{-A},\ \ \ A\, > 1\,;
\end{equation}
these are the numbers of irrationality measure greater than $2$. We shall treat the two cases separately.

We begin by re-writing the conclusion of corollary \ref{n=1} slightly. To simplify the appearance of various formulas, we shall suppose that $x\in\R$ satisfies
\begin{equation}\label{nongeneric1}
p/q\ > \ x\,;
\end{equation}
if $x>p/q$, the signs of various terms need to be changed. In any case, as $p/q$ runs through the continued fraction approximation of $x$, the sign of $p/q-x$ alternates. We also consider only the real part of $\phi_\tau$, without essential loss of generality. With these conventions,
\begin{equation}\label{nongeneric2}
\begin{aligned}
&\biggl| \Re\left( \phi_\tau(\textstyle\frac pq) \,-\,  \phi_\tau(x)\, + \, c_{\tau,0} \bigl(\textstyle\frac pq - x\bigr)\right)\, - \, \Re \biggl( \l^{-1}\,q^{-1}\,c_{\g\tau,0} (p- q x)^\l\biggr)
\\
&\qquad\qquad+\ \ \Re\biggl(\bigl(p - qx\bigr)^{\l+1}\phi_{\g\tau}(\g\,x)\biggr) \biggr|\ \leq \ C\,q \,\bigr(p - q x\bigl)^{\Re\l+2}\,.
\end{aligned}
\end{equation}
To get a handle on the behavior of $\phi_\tau$ near some irrational number $x_0$, we shall apply this estimate with $x=x_0$ fixed, and $p/q$ running through an approximating sequence of rational numbers. In this situation $\,\g$, $\phi_{\g\tau}$, and $c_{\g\tau,0}$ depend on the approximating rational $p/q$.

Thus let $p/q$ be an approximating sequence as in (\ref{generic}), converging to an irrational number $x_0$, such that $q>0$ and $p/q>x_0$. This latter assumption is merely a convenience; it will not affect our conclusions. We suppose that $\Re\phi_\tau$ satisfies a pointwise H\"older condition at $x_0$,
\begin{equation}\label{generic1}
|\,\Re\bigl(\phi_\tau(x)\,-\, \phi_\tau(x_0)\bigr)\,|\ \leq \ D\,|\,x\,-\, x_0 \,|^{\alpha}\,,
\end{equation}
for some $\alpha>0$. We shall see eventually that $\alpha=1$ is impossible. Following Duistermaat \cite{Du}, we pick a number $\eta$, $0<\eta<1$, whose value will depend on $q$, and define $x_\eta\in\R$ by the identity
\begin{equation}\label{generic2}
p\, - \,q\,x_\eta\ =\ \eta\,(\,p\,-\,q\,x_0).
\end{equation}
Then $x_\eta$ lies strictly between $x_0$ and $p/q$. At the very end of this section, we shall also work with $\eta>1$, which then reverses the order of $x_0$ and $x_\eta$.

At this point the argument branches. By going to a subsequence, we can arrange that either
\begin{equation}\label{genericbranch}
\begin{aligned}
&\text{a)\ \ $\Re c_{\g\tau,0}=0$\ \ for all terms of the sequence; or}\qquad\qquad\qquad\qquad
\\
&\text{b)\ \ $\Re c_{\g\tau,0}\neq 0$\ \ for all terms of the sequence}\,.
\end{aligned}
\end{equation}
In the situation a), we apply (\ref{generic1}) with $x=x_\eta$ and $x=p/q$, then add the two inequalities:
\begin{equation}\label{generic3}
\bigl|\Re\bigl(\phi_\tau(x_\eta)- \phi_\tau(p/q)\bigr)\bigr|\ \leq \ 2\,D\,|x_0- p/q|^{\alpha}\ = \ 2\,D\,q^{-\alpha}\,(p-qx_0)^\alpha\,.
\end{equation}
For future reference, we note that this remains correct for $\eta>1$, provided $x_0$ is replaced by $x_\eta$. We set $x=x_\eta$ in the inequality (\ref{nongeneric2}), and of course also $\Re c_{\g\tau,0}=0$, and combine the resulting inequality with (\ref{generic3}). This gives
\begin{equation}\label{generic4}
\begin{aligned}
\bigl| \Re\bigl( (p-qx_\eta)^{\l+1}\phi_{\g\tau}(\g\,x_\eta)\, + \, c_{\tau,0}\, q^{-1}(p - qx_\eta)\bigr)\bigr|\ \leq\ \frac{2\,D}{q^{\alpha}}\,(p-qx_0)^\alpha
\\
+ \ C\,q\,(p - qx_\eta)^{\Re\l+2}\,.
\end{aligned}
\end{equation}
Now divide by $(p-qx_\eta)^{\Re\l+1}$ and recall (\ref{generic}) and (\ref{generic2}):
\begin{equation}\label{generic5}
\begin{aligned}
\bigl| \Re\bigl( (p-qx_\eta)^{\Im\l}&\phi_{\g\tau}(\g\,x_\eta)\, + \, c_{\tau,0}\, q^{-1}(p - qx_\eta)^{-\Re\l}\bigr)\bigr|\ \leq
\\
&\leq\ 2\,D\,q^{-\alpha}\,\eta^{-\Re\l-1}\,(p-qx_0)^{\alpha-\Re\l-1} \ +\ C\,\eta\,.
\end{aligned}
\end{equation}
If $\Re \l\leq 0$, there exist choices of $\alpha$ such that $1+\Re \l \leq \alpha\leq 1$. In this situation (\ref{generic5}) will lead to a contradiction, thus ruling out (\ref{generic1}). Indeed, for $1+\Re \l \leq \alpha\leq 1$, (\ref{generic}) implies $(p-qx_0)^{\alpha-\Re\l-1}\leq q^{-\alpha+\Re\l+1}$, and (\ref{eitheror}) implies $c_{\tau,0}=0$, hence
\begin{equation}\label{generic6}
\begin{aligned}
\bigl| \Re\bigl( (p-qx_\eta)^{\Im\l}\phi_{\g\tau}(\g\,x_\eta)\bigr)\bigr|\ &\leq \ 2\,D\,q^{-2\alpha+\Re\l+1}\,\eta^{-\Re\l-1} \ +\ C\,\eta
\\
&\leq \ 2\,D\,q^{-\Re\l-1}\,\eta^{-\Re\l-1} \ +\ C\,\eta\,.
\end{aligned}
\end{equation}
We can bound the right hand side by any $M>0$, e.g.,
\begin{equation}\label{generic7}
M\ =\ \textstyle\frac 12\,\operatorname{min}\bigl\{ \operatorname{max}_{x\in\R}|\phi_{\g\tau}(x)|\, \mid\, \g\in SL(2,\Z)/\G\,\bigr\}\,,
\end{equation}
by restricting $\eta$ to the interval
\begin{equation}\label{generic8}
\bigl(\textstyle\frac{4D}{M}\bigr)^{\frac{1}{\Re \l +1}}\,q^{-1}\ \leq \ \eta\ \leq \ \operatorname{min}\left( \,\textstyle\frac 12\,,\, \frac M{2C}  \right)\,.
\end{equation}
As $x_\eta$ runs over the corresponding range, the argument of $(p-qx_\eta)^{\Im\l}$ covers an interval whose length grows like $\,\Im\l\,\, \log q\,$ as $q\to\infty$. Thus, if $\,\Im\l\neq 0$, one can argue as in the proof of corollary \ref{diffatrationals} and derive a contradiction from (\ref{generic1}); if $\,\Im\l= 0$ the argument becomes even simpler. In any case, this proves
\begin{equation}\label{generic9}
\begin{aligned}
&\text{if $\ \Re\,\tau \not\equiv 0$\,,\ \ $\Re\l \leq 0$\,,\ \ and $1 + \Re\l \leq \alpha$\,,}
\\
&\qquad{\lim\sup}_{x\to x_0}\,|x-x_0|^{-\alpha}\,|\Re\phi_\tau(x) - \Re\phi_\tau(x_0)|\ = \ \infty \,;
\end{aligned}
\end{equation}
in particular, $\Re\phi_\tau$ is not differentiable at $x=x_0$. The hypothesis $\,\Re\l\leq 0$ automatically puts us into the situation (\ref{genericbranch}a), so (\ref{generic9}) has been proved in complete generality.

We now use the fact that that $x_0$ has irrationality measure $2\,$: for any $\e>0$ and any approximating sequence $p/q$, one can arrange that
\begin{equation}\label{generic11}
q^{-1-\e}\ \leq \ |\,p\,-\,q\,x_0\,|\ < \ q^{-1}
\end{equation}
after dropping at most a finite number of terms. We allow positive values of $\l$, but do not exclude the case $\Re \l\leq 0$, which was already covered. If $\l>0$ we explicitly assume (\ref{genericbranch}a), which is otherwise automatic. We shall show that (\ref{generic1}) cannot hold if $\alpha > \frac{1+\Re\l}2$; this improves upon (\ref{generic9}) when $\Re\l\leq 0$ -- but only under the additional hypothesis (\ref{generic11}).

Indeed, if $\frac{1+\Re\l}2<\alpha<1+\Re\l$, we can choose some small $\e>0$, such that $(1+\e)(1+\Re\l-\alpha)-\alpha<0$ and $(1+\e)\Re\l-1<0$. That allows us to bound both $q^{-1}(p-qx_0)^{-\Re\l}$ and $q^{-\alpha}(p-qx_0)^{\alpha-\Re\l-1}$ by strictly negative powers of $q$, both for $\l>0$ and $\Re\l\leq 0$. The inequality (\ref{generic5}) now implies
\begin{equation}\label{generic12}
\bigl| \Re\bigl( (p-qx_\eta)^{\Im\l}\phi_{\g\tau}(\g\,x_\eta)\bigr) \bigr|\ \leq \ A\,q^{-a}\,\eta^{-\Re\l-1}\ + \ B\,q^{-b}\,\eta^{-\Re\l}\ + \ C\,\eta\,,
\end{equation}
with suitably chosen positive constants $A,\,B,\,C,\,a,\,b$. We can then argue as before: there exists an interval of choices for $x_\eta$, on which $\log(p-qx_\eta)$ runs through an interval of length proportional to $\log q$, and on which the right hand side of (\ref{generic12}) can be bounded by any given constant $M>0$. That is a contradiction unless $\,\tau\equiv 0$.

We continue with the hypothesis (\ref{generic11}), but now in the situation (\ref{genericbranch}b), which only occurs when $0<\l<1$. We can still argue as before, but we need to carry along $\l^{-1}q^{-1}\Re c_{\g\tau,0}(p-qx_\eta)^\l$ into the inequality (\ref{generic4}). After dividing by $(p-qx_\eta)^{1+\l}$, that leaves us with
\begin{equation}\label{generic13}
\begin{aligned}
\bigl| \Re\bigl(\phi_{\g\tau}(\g\,x_\eta) + \textstyle\frac  {c_{\tau,0}}{ q \eta^\lambda}&(p - qx_0)^{-\l} - \textstyle\frac{c_{\g\tau,0}}{\l q\eta}(p-qx_0)^{-1}\bigr)\bigr|\ \leq
\\
&\leq\ 2\,D\,q^{-\alpha}\,\eta^{-\l-1}\,(p-qx_0)^{\alpha-\l-1} \ +\ C\,\eta
\end{aligned}
\end{equation}
instead of (\ref{generic5}). Recall that the quantity $\,\Re
c_{\g\tau,0}$ takes on only finitely many values, all different
from zero by assumption. Hence
\begin{equation}\label{generic14}
\bigl|\l^{-1}\,\eta^{-1}\,q^{-1}\,(p-qx_0)^{-1}\, \Re c_{\g\tau,0}\bigr|\ > \ \eta^{-1}\,\l^{-1}\,\operatorname{min}_\g\,|\Re c_{\g\tau,0}|
\end{equation}
can be forced to tend to $\infty$ if we require $\eta<q^{-\e}$. From this point on we can argue essentially as before: fix $\a>\frac{1+\l}{2}$; if $\eta$ is bounded below by multiples of $q^{-\frac{a}{1+\l}}$ and of $q^{-\frac{b}{\l}}$~-- in the notation of (\ref{generic12})~-- and above by a multiple of $q^{-\e}$, the left hand side tends to $0$ as $p/q\to x_0$. That contradicts the boundedness of $\phi_{\g\tau}$, since the quantity (\ref{generic14}) is now unbounded. We have shown:
\begin{equation}\label{generic15}
\begin{aligned}
\text{if $\Re \tau \not\equiv 0$, if $\textstyle\frac{1+\Re\l}{2}<\a$, and if $x_0$ has irrationality measure $2$,}
\\
{\lim\sup}_{x\to x_0}\,|x-x_0|^{-\alpha}\,|\re\phi_\tau(x) - \Re\phi_\tau(x_0)|\ = \ \infty \,,
\end{aligned}
\end{equation}
for all possible choices of $\l$ and both cases of (\ref{genericbranch}).

We still need to treat the case (\ref{nongeneric}) of irrationality measure greater than $2$, but only for $0<\l<1$. We fix $x=x_0$ and choose an approxi\-mating sequence $p/q$, with $(p,q)=1$, $q>0$, subject to the condition (\ref{nongeneric}); as before, we also assume (\ref{nongeneric1}), without essential loss of generality. In the situation (\ref{genericbranch}b), we can work directly with (\ref{nongeneric2}) instead of (\ref{generic5}): we set $x=x_0$ in (\ref{nongeneric2}) and divide by $(p/q-x_0)^\alpha$, with $\l<\alpha\leq 1$. We then use the bound $|p-qx_0|<q^{-A}$ to estimate three of the resulting terms:
\begin{equation}\label{nongeneric3}
\begin{aligned}
&q^{\alpha-1} \bigl(p - q x_0\bigr)^{\l-\alpha}\, >\, q^{\alpha - 1 - A(\l-\alpha)}\, = \, q^{\alpha(A+1)-(A\l+1)}\,,
\\
&\qquad q^{\alpha} \bigl(p - q x_0\bigr)^{\l+1-\alpha}\, < \, q^{\alpha - A(\l+1-\alpha)}\, = \, q^{\alpha(A+1)-A(\l+1)}\,,
\\
&\qquad\qquad q^{1+\alpha}\bigl(p - q x_0\bigr)^{\l+2-\alpha}\, < \, q^{1 + \alpha - A(\l+2-\alpha)}\, = \, q^{\alpha(A+1)-A(\l+2)+1}.
\end{aligned}
\end{equation}
If $\alpha$ is picked so that $(A\l + 1)(A + 1)^{-1} < \alpha \ll 1$, the first of the expressions in (\ref{nongeneric3}) grows as $p/q\to x_0$, and the other two decay; also $\l<\alpha\leq 1$ in this situation, as was assumed. The quantity $\Re c_{\g\tau,0}$ can take on only finitely many different values. Thus, in view of (\ref{nongeneric2}) and (\ref{nongeneric3}),
\begin{equation}\label{nongeneric4}
\underset{p/q\to x_0}{\lim\operatorname{sup}}\, \bigl|\textstyle\frac pq -x_0\bigr|^{-\alpha}\, \bigl|\Re\bigl(\phi_\tau(\textstyle\frac pq) -  \phi_\tau(x_0)\bigr)\bigr|\, = \, \infty\ \ \ \text{if}\ \alpha>\frac{(A\l+1)}{A+1}\,,
\end{equation}
in the situation (\ref{genericbranch}b), with $A>1$ as in (\ref{nongeneric}). Since (\ref{nongeneric}) remains valid when $A$ decreases, we may let $A$ tend to $1$ from above. In other words, the conclusion of (\ref{nongeneric4}) applies for any $\alpha>\frac{\l+1}2$.

Still with $0<\l<1$, but now with $\,\Re c_{\g\tau,0}= 0$ for all terms of the approximating sequence, we return to the argument in (\rangeref{generic1}{generic8}), with $\eta>1$ instead of $0<\eta<1$ as before. We had remarked already that we need to replace $x_0$ by $x_\eta$ in (\ref{generic3}) and all corresponding quantities in the subsequent derivation. The analogue of (\ref{generic5}) in the current setting is
\begin{equation}\label{nongeneric5}
\begin{aligned}
&\bigl| \Re\phi_{\g\tau}(\g\,x_\eta)\bigr|\ \ \leq\ \ q^{-1}\,\eta^{-\l}\,(p - qx_0)^{-\l}\,| c_{\tau,0}|\, \ +
\\
&\qquad +\ 2\,D\,q^{-\alpha}\,\eta^{\alpha-\l-1}\,(p-qx_0)^{\alpha-\l-1} \ +\ C\,\eta\,q\,(p - qx_0)\,.
\end{aligned}
\end{equation}
To get this, we have left $q(p - qx_0)$ in its original form instead of bounding it by $1$, as in (\ref{generic5}). The additional $\alpha$ in the exponent of $\eta$ compensates for the change from $x_0$ to $x_\eta$ in (\rangeref{generic3}{generic4}), and we have moved the term involving $c_{\tau,0}$ to the other side of the inequality. We now define
\begin{equation}\label{nongeneric6}
\tilde\eta\ \ = \ \ \eta\,q\,(p-qx_0)\,.
\end{equation}
Since $q(p-qx_0)<q^{1-A}$, the hypothesis $\eta>1$ will be satisfied if $\tilde\eta>q^{1-A}$, as shall be assumed. With this new convention, the bound (\ref{nongeneric5}) becomes
\begin{equation}\label{nongeneric7}
\bigl| \Re\phi_{\g\tau}(\g\,x_\eta)\bigr|\ \ \leq\ \ | c_{\tau,0}|\,q^{\l-1}\,\tilde\eta^{-\l} \, + \, 2\,D\,q^{-2\alpha+\l+1}\,\tilde\eta^{\alpha-\l-1} \, +\, C\,\tilde\eta\,.
\end{equation}
We fix $\a$, $\frac{1+\l}2<\a\leq1$, and argue analogously to (\rangeref{generic7}{generic8}): there exists an interval of values $q^{-\e}\ll\tilde\eta\ll 1$ on which (\ref{nongeneric7}) is impossible; since
\begin{equation}\label{nongeneric8}
\textstyle\frac{d\ }{dx}\bigl( \g x \bigr)\bigr|_{x=x_\eta}\ = \ q^2\,\tilde\eta^{-2}\ \ \ \text{and}\ \ \ \textstyle\frac{d^2\ }{dx^2}\bigl( \g x \bigr)\,>0\,\ \ \text{for}\ \ x\,<\,p/q\,,
\end{equation}
$\g x_\eta$ runs over many cycles of the periodic function $\phi_\g$ as $\tilde\eta$ runs over its allowed range. That contradicts (\ref{nongeneric1}).

At this point we have established the non-differentiability of $\Re\phi_\tau$ at irrational points, in all possible cases. More precisely,

\begin{thm}\label{nondiffthm2}
If $x_0$ is irrational and $\,\Re\tau\not\equiv 0$,
\[
{\lim\operatorname{sup}}_{x\to x_0}\ |x - x_0|^{-\a}\,|\Re\phi_\tau(x) - \Re\phi_\tau(x_0)|\ = \ \infty
\]
in each of the following cases:\newline
\indent a)\ \ $\Re\l \leq 0$\ \ and\ \  $\a\geq 1+\Re\l$\,;\newline
\indent b)\ \ $\Re\l \leq 0$\,,\ \ $\a> \frac{1+\Re\l}2$\,,\ \ and\ \ $x_0$ has irrationality measure 2;\newline
\indent c)\ \ $\l > 0$\ \ and\ \ $\a> \frac{1+\Re\l}2$.
\end{thm}

\section{Final remarks}\label{final}

The group $G=SL(2,\R)$ has a twofold covering group $\widetilde G \to G$, the so-called metaplectic cover. The principal series of $\widetilde G$ is parameterized by pairs $(\l,\d)$ with $\l\in\C$ and $\d\in\Z/4\Z$, with $\d=\pm 1$ corresponding to ``genuine" representations of $\widetilde G$, i.e., representations that are not representations of $G$. The representations corresponding to $\d=0$ and $\d=2\,$ do drop down to representations of $G$, with $G$-parameters $(\l,\d/2)$ in the notation of $(\ref{vlambda})$. We can now consider $\widetilde\G$-automorphic distributions $\,\tau\in (V_{\l,\d}^{-\infty})^{\widetilde\G}$, for subgroups $\widetilde\G\subset\widetilde G$ whose image in $G$ is a normal subgroup of finite index $\G \subset SL(2,\Z)$. One does not get anything new except in the genuine case, and then only when $\widetilde\G$ does not contain the kernel of the covering morphism $\widetilde G \to G$.

When $\l=1/2$ and $\d=\pm 1$, automorphic distributions $\,\tau\in (V_{\l,\d}^{-\infty})^{\widetilde\G}$ correspond to modular forms of weight $1/2$ -- on the upper half plane for one of the two choices of $\d$, and on the lower half plane for the other $\d$. A modular form of this type is square integrable if and only if $\,\tau$ satisfies the condition analogous to (\ref{eitheror}); we shall assume this is the case. One can then identify $\,\tau$ with a distribution on $\,\R$, as before. Since $\widetilde \G\cong \G$ via the covering morphism $\widetilde G \to G$, one can even characterize the $\widetilde\G$-invariance in terms of the action of $\G$ on $\,\R\cup \{\infty\}$\,: when $\g^{-1}\in\G$ has matrix entries $a,b,c,d$ as in (\ref{pilambda}), invariance under the lifted $\widetilde\g\in \widetilde \G$ translates into the identity
\begin{equation}\label{vlambdatilde}
\begin{aligned}
\tau(x)\ = \ \chi\,&|cx+d|^{1/2}\,\tau\left( \frac{ax+b}{cx+d} \right),
\\
&\text{with}\ \ \chi = \chi\bigl(\g,\d,\sg(cx+d)\bigr)\in\{\pm 1, \pm i\}\,;
\end{aligned}
\end{equation}
because of the hypothesis (\ref{eitheror}), this has a definite meaning even at $x=-c/d$. Iwaniec \cite[\S 2.8]{Iw} treats modular forms of weight 1/2. The identity (\ref{vlambdatilde}) follows from his discussion by taking the distribution limit as $y\to 0$.

The theta series $\theta(z)= \sum_{n\in\Z} e(n^2z)$ is the most important modular form of weight $1/2$. Like $\theta(z)$ itself, the boundary distribution
\begin{equation}\label{theta1}
\theta(x)\ \ = \ \ {\sum}_{n\in\Z}\ e(n^2x)
\end{equation}
is invariant under $\widetilde\G_1(4)$ and transforms under $\widetilde\G_0(4)$ according to the non-trivial character of $\widetilde\G_0(4)/\widetilde\G_1(4)\cong \{\pm 1\}$. Let
\begin{equation}\label{theta2}
\phi_\theta(x)\ \ = \ \ \frac 1{2\pi i}\ {\sum}_{n\neq 0}\ n^{-2}\,e(n^2x)
\end{equation}
denote the antiderivative of $\theta(x)-1$, in analogy to (\ref{contantider1}). Then $\pi \Re\phi_\theta(x)$ coincides with the function (\ref{Riemann}). In view of the invariance conditions (\ref{vlambdatilde}), the arguments of sections \ref{ratsec} and \ref{irratsec} apply without change; in this situation, they become essentially equivalent to Duistermaat's argument \cite{Du}, of course. The group $\G_0(4)$ has exactly three orbits in $\Q\cup\{\infty\}$. The numbers $p/2q$, with $p,q$ odd, constitute one orbit. At the points of the other two orbits, the non-cuspidal nature of $\theta(x)$ prevents $ \Re\phi_\theta(x)$ from being differentiable.

The infinite product $\eta(z)=e(z/24)\prod_{n=1}^\infty (1 - e(nz))$ is another modular form of weight $1/2$, modular with respect to $\widetilde\G_1(24)$, and transforming according to a character of order $24$ under the action of $\widetilde{SL}(2,\Z)$. Unlike $\theta(z)$, it is cuspidal. The antiderivative $\phi_\eta(x)$ of the boundary distribution $\eta(x)$ therefore has derivative zero at all rational points. According to theorem \ref{nondiffthm2}, the real and imaginary parts of both $\phi_\theta$ and $\phi_\eta$ violate the pointwise H\"older condition with any index $\a>3/4$, at any irrational point. In the case of $\phi_\eta$, this statement is absolutely sharp, since theorem \ref{regthm} asserts global H\"older continuity of index $3/4\,$; $\,\phi_\theta$, on the other hand, is only H\"older continuous of index $1/2$.

\begin{figure}
\begin{center}
\includegraphics[height=2in,width=5.3in]{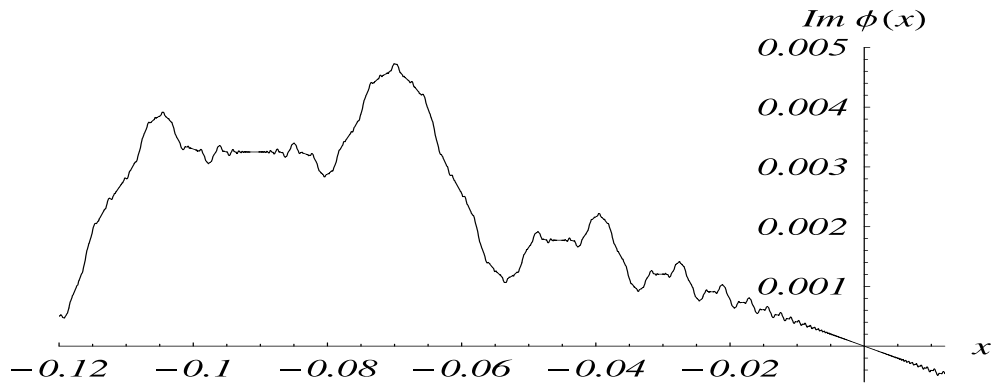}
\end{center}
\caption{The imaginary part of the second antiderivative of the automorphic distribution corresponding to the modular form $F(z)= \eta(z)^2 \eta(11z)^2$.
}\label{imagweight2medium+_3}
\end{figure}

The function $F(z)= \eta(z)^2 \eta(11z)^2$ is a cuspidal modular form of weight two for $\G=\G_0(11)$. Alternatively $F(z)$ can be described as the modular form associated to the elliptic curve
\begin{equation}\label{ellcu2}
y^2 \ + \ y \ \ = \ \ x^3 \ - \ x^2 \ - \ 10 \, x \ - \ 20\,.
\end{equation}
As in the case of all modular forms of weight two, one needs to take two antiderivatives of the the boundary distribution $F(x)$ to get a continuous function. That function is differentiable at rational points, but both its real and imaginary part are non-differentiable at irrational points, as can be shown by an adaptation of the arguments in section \ref{irratsec}. Figure \ref{imagweight2medium+_3} shows the imaginary part.
\bigskip
\bigskip

\begin{tabular}{lcl}
Stephen D. Miller                    & & Wilfried Schmid \\
Department of Mathematics            & & Department of Mathematics \\
Hill Center-Busch Campus             & & Harvard University \\
Rutgers University                        & & Cambridge, MA 02138 \\
110 Frelinghuysen Rd                 & & {\tt schmid@math.harvard.edu}\\
Piscataway, NJ 08854-8019            & & \\
{\tt miller@math.rutgers.edu}
\end{tabular}

\end{document}